\documentstyle[11pt,amscd]{amsart}
\newtheorem{thm}{Theorem}[section]
\newtheorem{lem}[thm]{Lemma}
\newtheorem{cor}[thm]{Corollary}
\newtheorem{prop}[thm]{Proposition}

\theoremstyle{definition}

\def \CPb {\overline{\bold{CP}}^{\,2}}
\def \z {\zeta}
\def \R {\bold{R}}
\def \Z {\bold{Z}}

\def \CM {\cal M}

\def \asd {anti-self-dual }

\def \la {\langle}
\def \ra {\rangle}
\def \a {\alpha}
\def \b {\beta}

\def \lam {\lambda}

\def \o {\omega}
\def \s {\sigma}
\def \t {\tau}

\def \bd {\partial}

\def \x {\times}
\def \ve {\varepsilon}
\def \D {\bold{D}}
\def \A {\bold{A}}

\begin{document}

\baselineskip.5cm

\title{The Blowup Formula for Donaldson Invariants}
\author[Ronald Fintushel]{Ronald Fintushel}
\address{Department of Mathematics, Michigan State University \newline
\hspace*{.375in}East Lansing, Michigan 48824}   \email{ronfint@@math.msu.edu}
\thanks{The first author was partially supported NSF Grant DMS9102522  and the
second author by NSF Grant DMS9302526}  \author[Ronald J. Stern]{Ronald J.
Stern}
\address{Department of Mathematics, University of California \newline
\hspace*{.375in}Irvine,  California 92717}  \email{rstern@@math.uci.edu}
\date{June 14, 1994}

\thanks{The authors wish to thank Peter Kronheimer and Tom Mrowka for many
useful
suggestions and Tom Parker for helpful conversations. We also thank those,
especially Thomas Leness
and Paolo Lisca, who have sent us comments concerning an earlier version of
this paper.}

\maketitle

\section{Introduction}

Since their introduction in 1984 \cite{Donpoly}, the Donaldson invariants of
smooth $4$-manifolds
have remained as mysterious as they have been useful. However, in the past year
there has been a surge
of activity pointed at comprehension of the structure of these invariants
\cite{KM,FS}. One key to
these advances and to future insights lies in understanding the relation of the
Donaldson invariants
of a $4$-manifold $X$ and those of its blowup $\hat{X}=X\#\CPb$. It is the
purpose of this paper to
present such a blowup formula. This formula is independent of $X$ and is given
in terms of an
infinite series
\[ B(x,t) = \sum\limits_{k=0}^{\infty}B_k(x){t^k\over k!} \]
which is calculated in \S4 below.
This formula has been the target of much recent work.
The abstract fact that there exists such a formula which is
independent of $X$ was first proved by C. Taubes using techniques of
\cite{Reds}. J. Bryan \cite{B}
and P. Ozsvath \cite{Ozs} have independently calculated the coefficients
through $B_{10}(x)$. Quite
recently, J. Morgan and Ozsvath have announced a scheme which can recursively
compute all of the
$B_k(x)$. The special case of the blowup formula for manifolds of ``simple
type'' (see \S5 below) was
first given by P. Kronheimer and T. Mrowka. However, none of the techniques in
these cases approach
the simplicity of that offered here.

Before presenting the formula, we shall first establish notation for the
Donaldson invariants of a
simply connected $4$-manifold $X$ with $b^+>1$ and odd. (The hypothesis of
simple
connectivity is not necessary, but makes the exposition easier.) An orientation
of $X$, together with
an orientation of $H^2_+(X;\R)$ is called a {\em homology orientation} of $X$.
Such a homology
orientation determines the ($SU(2)$) Donaldson invariant, a linear function
\[ D=D_X:\A(X)=\text{Sym}_*(H_0(X)\oplus H_2(X))\to \R \]
which is a homology orientation-preserving diffeomorphism invariant.
Here
\[\A(X)=\text{Sym}_*(H_0(X)\oplus H_2(X))\]
is viewed as a graded algebra where $H_i(X)$ has degree $\frac12(4-i)$.
We let $x\in H_0(X)$ be the generator $[1]$ corresponding to
the orientation. Then as usual, if $a+2b=d>\frac34(1+b^+_X)$ and $\a\in
H_2(X)$,
\[ D(\a^ax^b)=\la\mu(\a)^a\nu^b,[\CM_X^{2d}]\ra \]
where $[\CM_X^{2d}]$ is the fundamental class of the (compactified)
$2d$-dimensional moduli space of
\asd connections on an $SU(2)$ bundle over $X$, $\mu:H_i(X)\to H^{4-i}(\cal
B^*_X)$ is the
canonical map to the cohomology of the space of irreducible connections on that
bundle
\cite{Donpoly}, and $\nu=\mu(x)$. The extension of the definition to smaller
$d$ is given in
\cite{MMblowup} (and is accomplished, in fact, from the knowledge of the lowest
coefficient in the
$SO(3)$ blowup formula). Since an $SU(2)$ bundle $P$ over $X$ has a moduli
space of dimension
\[\dim\CM_X(P) = 8c_2(P)-3(1+b^+_X) \]
it follows that such moduli spaces $\CM_X^{2d}$ can exist only for
$d\equiv\frac12(1+b^+_X)\pmod4$. Thus the Donaldson invariant $D$ is defined
only on elements of
$\A(X)$ whose total degree is congruent to $\frac12(1+b^+_X)\pmod4$.
By definition, $D$ is $0$ on all elements of other degrees.

We can now state the blowup formula. Let $\hat{X}=X\#\CPb$ and let $e\in
H_2(X)$ denote the homology
class of the exceptional divisor. Since $b^+_X=b^+_{\hat{X}}$, the
corresponding Donaldson invariants $D=D_X$ and $\hat{D}=D_{\hat{X}}$ have their
(possible) nonzero
values in the same degrees $\pmod4$. We first show that there are polynomials
$B_k(x)$ satisfying
\[ \hat{D}(e^k\,z) = D(B_k(x)\,z) \]
for all $\,z\in \A(X)$ and then define the formal power series $B(x,t)$ as
above. Our result is that
\[ B(x,t)=e^{-{t^2x\over6}}\s_3(t) \]
where $\s_3$ is a particular quasi-periodic Weierstrass
sigma-function \cite{Ak} associated to the $\wp$-function which
satisfies the differential equation
\[(y')^2=4y^3-g_2y-g_3 \]
where
\[ g_2=4\,(\frac{x^2}{3}-1)\, , \hspace{.25in} g_3={8x^3-36x\over 27} \,. \]

There are also Donaldson invariants associated to $SO(3)$ bundles $V$ over $X$.
To define these
invariants one needs, along with a homology orientation of $X$, an integral
lift of $w_2(V)$. If
$c\in H_2(X;\Z)$ is the Poincar\'e dual of the lift, the invariant is denoted
$D_c$ or $D_{X,c}$ if
the manifold $X$ is in doubt. $D_c$ is nonzero only in degrees congruent to
$-c\cdot c+\frac12(1+b^+)\pmod4$. If $c'\equiv c\pmod2$ then
\[ D_{c'}=(-1)^{({c'-c\over2})^2}D_c\,. \]

The $SO(3)$ blowup formula  states that there are polynomials $S_k(x)$ such
that
\[ \hat{D}_e(e^k) = D(S_k(x)) \]
and if
\[ S(x,t) = \sum\limits_{k=0}^{\infty}S_k(x){t^k\over k!} \]
then
\[ S(x,t)=e^{-{t^2x\over6}}\s(t) \] where $\s(t)$  is the standard
Weierstrass
sigma-function \cite{Ak} associated to $\wp$.
The coefficients $S_k(x)$  for $k\le 7$ were earlier computed by T. Leness
\cite{Leness}.

The discriminant of the cubic equation
$4y^3-g_2y-g_3=0$ turns out to be ${x^2-4\over 4}$. Thus, when (viewed as a
function on
$\A(X)$) $D(x^2-4)=0$, the Weierstrass sigma-functions degenerate to elementary
functions, and the
blowup formula can be restated in terms of these functions. This is done in the
final section. It is
interesting to note that the condition $D(x^2z)=4\,D(z)$ is the {\em simple
type} condition introduced
by Kronheimer and Mrowka \cite{KM}.

Our formulas are proved by means of a simple relation satisfied by  $D(\t^4z)$
where $\t\in
H_2(X;\Z)$ is represented by an embedded $2$-sphere of self-intersection $-2$.
When this
relation is applied to $\t=e_1-e_2$, the difference of the two exceptional
classes of the double
blowup $X\#2\CPb$, one obtains a differential equation for $B(x,t)$. Solving
this equation gives our
formulas.

\section{Some Relations among Donaldson Invariants}

The key to the blowup formula lies in a few simple  relations which are useful
for evaluating
Donaldson invariants on classes represented by embedded spheres of
self-intersection $-2$ and $-3$.
We begin by studying the behavior of the Donaldson invariant of a  $4$-manifold
with a homology class $\t$ represented by an embedded $2$-sphere $S$ of
self-intersection
$\t\cdot \t=-2$.  Let $\la \t\ra^\perp$ denote $\{\a\in H_2(X)|\t\cdot\a=0\}$
and let
\[\A(\t^\perp)=\A_X(\t^\perp)=\text{Sym}_*(H_0(X)\oplus \la \t\ra^\perp)\,.\]

\begin{thm}{\em (Ruberman\; \cite{R})} Suppose that $\t\in H_2(X;\Z)$ with
$\t\cdot \t=-2$ is represented by an embedded sphere $S$. Then for $\,z \in
\A(\t^\perp)$, we have
$D(\t^2\,z)=2\,D_{\t}(\,z)$. \label{Ruber}\end{thm}

\begin{pf} Write $X=X_0 \cup N$ where $N$ is a tubular
neighborhood of $S$, and note that $\bd N$ is the lens space $L(2,-1)$. Since
$b^+_{X_0}>0$, generically there are no reducible \asd connections on $X_0$.
However, since $b^+_{N}=0$, there are nontrivial reducible \asd connections
arising
from complex line bundles  $\lam^m$, $m \in\Z$, where $\la c_1(\lam),\t\ra=-1$.
The corresponding moduli spaces $\CM_N(\lam^m\oplus\bar{\lam}^m)$ have
dimensions
$4m^2-3$ and have boundary values $\z^m$ where $\z$ generates the character
variety of
$SU(2)$ representations of $\pi_1(\bd N)=\Z_2$ mod conjugacy. (Of course,
$\z^{2m}$
is trivial, and $\z^{2m+1}=\z$.)

Since $\la \t\ra^\perp=H_2(X_0)$, we need to evaluate the Donaldson
invariant on two copies of $\t$ and classes in $H_2(X_0)$. After cutting down
moduli spaces by
intersecting with transverse divisors representing the images under $\mu$ of
these classes in
$H_2(X_0)$ and using the given homology orientation, we may assume without loss
of generality that
there are no such classes and that we are working with a  $4$-dimensional
moduli space $\CM_X$. Let
$V_1$ and $V_2$ be divisors representing $\mu(\t)$, coming from general
positioned surfaces in $N$.
The Donaldson invariant is the signed intersection number
\[ D(\t^2)=\#(\CM_X\cap V_1\cap V_2). \]
A standard dimension counting argument (cf. \cite{Donpoly}) shows that if we
choose a metric on $X$
with long enough neck length, $\bd N\x [0,T]$, then
all the intersections take place in  a neighborhood $\cal{U}$ of the grafted
moduli
space   $\CM_{X_0}[\z]\# \{A_{\lam}\}$   where $A_{\lam}$ is the reducible \asd
connection on
$\lam\oplus \bar{\lam}$, and $\CM_{X_0}[\z]$ is the $0$-dimensional cylindrical
end moduli space on
$X_0$ consisting of \asd connections which decay exponentially to the boundary
value $\z$. Let
$m_{X_0}$ be the signed count of points in $\CM_{X_0}[\z]$. A  neighborhood of
$A_{\lam}$ in the
moduli space $\CM_N(\lam\oplus\bar{\lam})$ is diffeomorphic to  $(\bold C
\x_{S^1} SO(3))/SO(3)
\cong\bold C/S^1\cong[0,\infty)$. Here $S^1$ acts on $SO(3)$ so that
$SO(3)/S^1=S^2$ and on $\bold C$
with weight $-2$. Thus the neighborhood $\cal{U}$ is
\[(\tilde{\CM}_{X_0}[\z]\x (\bold C \times_{S^1} SO(3)))/SO(3)\] where
``$\tilde{\CM}_{X_0}[\z]$''
denotes the based moduli space.

Now $\tilde{V}_1\cap (\bold C \times_{S^1} SO(3))=\{0\} \times_{S^1} SO(3)$,
and the
intersection of $V_1$ with all of $\CM_X$ is
\[(\CM_{X_0}[\z] \x(\{0\} \times_{S^1} SO(3)))/SO(3)=\Delta .\]
Fix a point $p\in\CM_{X_0}[\z]$, let $SO(3)\cdot p$ denote its orbit in
$\tilde{\CM}_{X_0}[\z]$,
and let
\[\Delta_p=SO(3)\cdot p\x (\{0\} \times_{S^1} SO(3)))/SO(3)\cong S^2.\]
Identify $\Delta_p$ with a transversal in $\tilde{\Delta}_p$ and compute the
intersection number
$\tilde{V}_2\cdot \Delta_p=\iota_p$. Since $\iota_p$ is independent of $p\in
\CM_{X_0}[\z]$, we
have  $D(\t^2)=\iota_p\cdot m_{X_0}$. The constant $\iota_p$ is
computed in \cite{FMbook} as follows. Note that $\Delta_p=\{0\}\times_{S^1}
SO(3)
\subset  \bold C\times_{S^1} SO(3)$ is a \,zero-section of the $c_1=-2$ complex
line
bundle over $S^2$ and $\tilde{V}_2$ is another section. Thus
$\tilde{V}_2\cdot \Delta_p=-2$; and so  $D(\t^2)=-2\,m_{X_0}$.

To identify the relative invariant $m_{X_0}$, view $\CM_{X_0}[\z]$ as
$\CM_{X_0,0}[\text{ad}(\z)]$, an $SO(3)$ moduli space. Since $\text{ad}(\z)$ is
the
trivial $SO(3)$-representation, we may graft connections in
$\CM_{X_0,0}[\text{ad}(\z)]$ to the trivial $SO(3)$ connection over $N$, and
since
$b_N^+=0$, there is no obstruction to doing this. We obtain an $SO(3)$ moduli
space
over $X$ corresponding to an $SO(3)$ bundle over $X$ with $w_2$
Poincar\'e dual to $\t$. (This is the unique nonzero class in $H^2(X;\Z_2)$
which
restricts trivially to both $N$ and $X_0$.)  Thus for   $\,z \in \A(\t^\perp)$,
we have
$D(\t^2\,z)=\pm2\,D_{\t}(\,z)$. (Note that since $\t\cdot\t=-2$, we have
$D_{-\t}=D_\t$.)

To determine the sign in this equation, we need to compare orientations on the
moduli spaces which
are involved. Let $A_0\in \CM_{X_0}[\z]$. The way that a sign is attached
to this point is described in \cite{Donor,K}. By addition and subtraction of
instantons, $A_0$ is
related to a connection $B_0$ in a reducible bundle $E$ over $X_0$, and $B_0$
can be connected by a
path to a reducible connection $R$ which comes from a splitting $E\cong
L_0\oplus\bar{L}_0$.
There is a standard orientation for the determinant line of
the operator $d_R^+\oplus d_R^*$, and this can be followed back to give an
orientation
for the determinant line at $A_0$. This determinant line is canonically
oriented because the
cohomology $H^*_{A_0}$ vanishes. Comparing the two orientations gives a sign,
$\ve$.

To determine the sign at the grafted connection $A_0\# A_{\lam}$, note that the
same sequence of
instanton additions and subtractions as above relates  $A_0\# A_{\lam}$ to
$B_0\# A_{\lam}$ which can
be connected to $R\# A_{\lam}$, a reducible connection on the bundle
$L\oplus\R$ over $X$, where
the Mayer-Vietoris map $H^2(X)\to H^2(X_0)\oplus H^2(N)$ carries $c=c_1(L)$ to
$c_1(L_0)+c_1(\lam)$.
Since $R\# A_{\lam}$ is reducible, there is an orientation of the determinant
line, and it relates to the orientation which can be pulled back from the
trivial connection by
$(-1)^{c\cdot c}$. Thus pulling the orientation back over $A_0\# A_{\lam}$
gives the
sign $\ve\cdot (-1)^{c\cdot c}$.

To get the sign for $A_0\#\Theta$ we first pass to $SO(3)$, and then
$\text{ad}(A_0)$ is related as
above to the reducible connection $\text{ad}(R)$ which lives in the line bundle
$L_0^2$.
Grafting to the trivial connection $\Theta_N$, we get
$\text{ad}(A_0)\#\Theta_N$ which is
connected to the reducible connection $\text{ad}(R)\#\Theta_N$. This lives in
the grafted line bundle
$L_0^2\#\R$ which has $c_1=2c_1(L_0)$. (Note that although $c_1(L_0)$ is not a
global class,
$2c_1(L_0)$ is.) The class $2c_1(L_0)$ restricts trivially to $X_0$ and to $N$
(mod $2$); so its
mod $2$ reduction is the same as that of $\t$. (We are here identifying $\t$
and its Poincar\'e
dual.) Since $\t=2c_1(\lam)$, the difference in these reductions is
$2c_1(L_0)-\t=2(c-\t)$. The
corresponding orientations compare via the parity of $(c-\t)\cdot(c-\t)\equiv
c\cdot c\pmod2$. Thus
the sign which is attached to  $A_0\#\Theta$ is $\ve\cdot (-1)^{c\cdot c}$, the
same as for $A_0\#
A_{\lam}$, and the sign in the formula above is `${}\,+\,{}$'. \end{pf}

For the case of the $SO(3)$ invariants the proof of Theorem~\ref{Ruber} can be
easily adapted to show:

\begin{thm} Suppose that $\t\in H_2(X;\Z)$ with $\t\cdot\t=-2$ is represented
by an
embedded sphere $S$.  Let $c\in H_2(X;\Z)$ satisfy $c\cdot\t\equiv0\pmod2$.
Then for
 $\,z \in \A(\t^\perp)$ we have
$D_c(\t^2\,z)=2\,D_{c+\t}(z)$. \ \ \ \qed
\label{RuberSO3}\end{thm}

We next need to review some elementary facts concerning the Donaldson
invariants of blowups. These
can be found, for example in \cite{FMbook,Ko,Leness}. Let $X$ have the
Donaldson invariant $D$, and
let $\hat{X}=X\#\CPb$ have the invariant $\hat{D}$.

\begin{lem} Let $e\in H_2(\CPb;\Z)\subset H_2(\hat{X};\Z)$ be the exceptional
class, and let
$c\in H_2(X;\Z)$. Then for all $\,z\in \A(X)$:
\begin{enumerate}

  \item $\hat{D}_c(e^{2k+1}\,z)=0$ for all $k\ge 0$.
  \item $\hat{D}_c(\,z)=D_c(\,z)$.
  \item $\hat{D}_c(e^2\,z)=0$.
  \item $\hat{D}_c(e^4\,z)=-2\,D_c(\,z)$.
	\item $\hat{D}_{c+e}(e^{2k}\,z)=0$ for all $k\ge 0$.
  \item $\hat{D}_{c+e}(e\,z)=D_c(\,z)$.
  \item $\hat{D}_{c+e}(e^3\,z)=-D_c(x\,z)$.
 \end{enumerate} \label{blowuplow} \end{lem}

\begin{pf} Items (1)--(5) are standard and are explained in \cite{FMbook}. Both
(1) and (5)
follow because the automorphism of $H_2(X\#\CPb;\Z)$ given by reflection in $e$
is realized by a
diffeomorphism.  Items (2) and (3) follow from counting arguments, and (4)
follows
from simple arguments as in the proof of Theorem~\ref{Ruber} above
\cite{FMbook}. Item (6) is due to
D. Kotschick \cite{Ko}.

A proof of (7) is given in \cite{Leness}. (However, the sign there differs from
ours since item
(6) is stated in \cite{Leness} with an incorrect sign.) We sketch a proof here.
Consider a
neighborhood $N$ of the exceptional curve, and let $X_0=\hat{X}\setminus N
\cong X\setminus B^4$. As
in the proof of Ruberman's theorem we lose no generality by assuming that we
are evaluating $\hat{D}$
only on $e^3$. A dimension counting argument shows that if we stretch the neck
between $X_0$ and $N$
to have infinite length by taking a sequence of generic metrics $\{ g_n\}$, and
if $V_i$ are
transverse divisors representing $\mu(e)$, then any sequence of connections  \[
A_n\in
\CM_{\hat{X},c+e}(g_n)\cap V_1\cap V_2\cap V_3\]  must converge to the sum of a
connection in the
$4$-manifold $\CM_{X_0,c}$ and the unique reducible connection on $N$
corresponding to the line bundle $\lam$ over $N$ whose Euler class is
Poincar\'e dual to $e$. The
based moduli space $\tilde{\CM}_N(\lam)$ is the orbit of this reducible
connection, a $2$-sphere,
$S^2_e$. Let $v$ denote the (positive) generator of the equivariant cohomology
$H^*_{SO(3)}(S^2_e)\cong H^*(\bold{CP}^{\infty})$ in dimension $2$. The class
$\mu(e)$ lifts to the
equivariant class  $-\frac12\la c_1(\lam),e\ra\,v=\frac12 v\in
H^2_{SO(3)}(S^2_e)$. The connections
in $S^2_e$ are asymptotically trivial and this induces an $SO(3)$ equivariant
push-forward map
\[ \bd_*(N): H^*_{SO(3)}(S^2_e)\to H^*_{SO(3)}({1})=H^*(BSO(3)).\]
If $u\in H^*_{SO(3)}({1})$ is the generator in dimension $4$ then
$\bd_*(N)(v^{2k+1})=2\,u^k$. So $\bd_*(N)(e^3)=\frac14 u$.

Since each connection in $\CM_{X_0,c}$ is also asymptotically trivial, there is
an induced map
$\bd^*(X_0):H^*_{SO(3)}({1})\to H^*_{SO(3)}(\tilde{\CM}_{X_0,c})$.
It follows from \cite{Reds,AB} that $\hat{D}_{c+e}(e^3)$ is obtained by
evaluating
\[\la \bd^*(X_0)\,\bd_*(N)(v^3),[\tilde{\CM}_{X_0,c}]\ra =
 \frac14\,\la \bd^*(X_0)(u),[\tilde{\CM}_{X_0,c}]\ra
 = \frac14\,\la \pi_*\bd^*(X_0)(u), [\CM_{X_0,c}] \ra \]
where basepoint fibration $\b$ over $X_0$ is
\[\pi:\tilde{\CM}_{X_0,c}\to\CM_{X_0,c}\, ,\]
the last equality because the $SO(3)$ action on $\tilde{\CM}_{X_0,c}$ is free.
But Austin and Braam \cite{AB}, for example, show that
$\pi_*\bd^*(X_0)(u)=p_1(\b)$. Since
$\nu=-\frac14 p_1(\b)$, we get  $\la\nu,[\CM_{X_0,c}]\ra = -D_c(x)$.
\end{pf}

We next consider embedded $2$-spheres of self-intersection $-3$.

\begin{thm} Suppose that $\t\in H_2(X;\Z)$ is represented by an
embedded 2-sphere $S$ with self-intersection $-3$. Let $\o\in H_2(X;\Z)$
satisfy
$\o\cdot\t\equiv0\pmod2$. Then for all $\,z\in\A(\t^\perp)$ we have
\[D_\o(\t\,z) = -D_{\o+\t}(\,z).\]
\label{3curve}\end{thm}
\begin{pf} The proof is similar in structure to that of
Theorem~\ref{Ruber}. Write $X=X_0\cup N$ where $N$ is a tubular neighborhood of
$S$.
Then $\bd N=L(3,-1)$. Let $\eta$ generate the character variety of $SO(3)$
representations of $\pi_1(\bd N)$. Reducible \asd $SO(3)$ connections on $N$
arise
from complex line bundles $\lam^m$, $m\in\Z$, where $\la c_1(\lam),\t\ra=-1$.
The corresponding moduli spaces $\CM_N(\lam^m\oplus\R)$ have boundary values
$\eta^m$
and dimensions $\frac{2}{3}m^2-3$ if $m\equiv 0\pmod3$ and ${2m^2+1\over3}-2$
if
$m\not\equiv0\pmod3$.

Since it is easiest to work with an $\o$ which satisfies $\bd\o_{X_0}=0\in
H_1(\bd X;\Z)=\Z_3$, we simply work with $\rho=3\,\o$ rather than $\o$. This is
no
problem, since $D_{3\o}=(-1)^{\o\cdot\o}D_{\o}$. Thus we may write
$\rho=\rho_0+\rho_N\in H_2(X_0;\Z)\oplus H_2(N;\Z)$.
As in our previous arguments, we assume that we are evaluating $D_{\rho}$ only
on $\t$.
A dimension counting argument
shows that $D_{\rho}(\t)$ is the product of relative invariants
$D_{X_0}[\eta^m]$
coming from a $0$-dimensional cylindrical end moduli space over $X_0$
with terms coming from nontrivial reducible
connections on $N$. These reducible connections must live in moduli spaces of
dimension $\le 1$, and the corresponding line bundles must have $c_1\equiv
\rho_N\pmod2$. Our hypothesis, $\o\cdot\t\equiv0\pmod2$ implies that
$\rho_N\cdot\t\equiv0\pmod2$; so the line bundle in question must be an even
power of
$\lam$. Recalling the constraint that the dimension of the corresponding moduli
space be $\le1$, the only possibility is $\CM_N(\lam^2\oplus\R)$.

Consider an \asd connection $A_0$ lying in the finite $0$-dimensional moduli
space $\CM_{X_0}[\eta^2]$, and let $A_{\lam^2}$ be the reducible \asd
connection on $N$. A neighborhood of the $SO(3)$ orbit of $A_{\lam^2}$ in the
based
moduli space   $\tilde{\CM}_N(\lam^2\oplus\R)$ is modelled by
$SO(3)\x_{S^1}\bold{C}$ and the (based) divisor for  $\t$ is
$-\frac12 \la c_1(\lam^2),\t\ra (SO(3)\x_{S^1}\{0\})=SO(3)\x_{S^1}\{0\}$.
The based connections obtained from grafting the orbit $SO(3)_{A_0}$ of $A_0$
to the orbit of
$A_{\lam^2}$ are given by the fibered product of these orbits over the 2-sphere
in
$SO(3)$ consisting of representations of $\pi_1(\bd N)$ which are in the
conjugacy
class $\eta^2$. By cutting this down by the divisor for $\t$ we obtain (up to
sign) the fibered product of $SO(3)_{A_0}$ with $S^2$ over $S^2$; i.e. simply
$SO(3)_{A_0}$. Taking the quotient by $SO(3)$,
\[ D_{\rho}(\t)=\pm D_{X_0}[\eta^2].\]

Since $\eta^2=\eta$ in the character variety (a copy of $\Z_3$), we can
graft \asd connections $A_0$ to the unique (reducible) connection $A_{\lam}$
lying in the
moduli space $\CM_N(\lam\oplus\R)$ of formal dimension $-1$. As the
glued-together
bundle has $w_2$ which is Poincar\'e dual to $\rho+\t\pmod2$, we have
\[ D_{\rho_0}[\eta^2]=\pm D_{\rho+\t},\]
so our result is proved up to a sign.

To get this sign, we need to compare signs induced at $A_0\# A_{\lam^2}$ and
$A_0\# A_{\lam}$ using
a fixed homology orientation of $X$ and the integral lifts $\rho$ and $\rho+\t$
of the
corresponding Stiefel-Whitney classes. By an excision argument \cite{Donor},
the difference in signs
depends only on the part of the connections over the neighborhood $N$. Thus the
sign is universal,
and may be determined by an example. For this, let $X$ be the $K3$ surface and
$\hat{X}=X\#\CPb$.
Let $s$ be any class in $H_2(X)$ of square $-2$ represented by an embedded
$2$-sphere (e.g. a
section), and let $\t=s+e$. Note that $s-2e\in\A(\t^{\perp})$.
Then using Theorems \ref{Ruber} and \ref{blowuplow},
\begin{eqnarray*}
\hat{D}((s-2e)\t)&=&D(s^2)=2\,D_s\\
\hat{D}_{\t}(s-2e)&=&-2\,D_s
\end{eqnarray*}
so the overall sign is `${}\,-\,{}$'.
\end{pf}

Next we combine our two relations to obtain a relation which is crucial in
obtaining
the general blowup formula. This relation was first proved by Wojciech
Wieczorek
using different methods. His proof will appear in his thesis \cite {Wiecz}.

\begin{cor} Suppose that $\t\in H_2(X;\Z)$ is represented by an embedded
$2$-sphere with
self-intersection $-2$, and let $c\in H_2(X;\Z)$ with $c\cdot\t\equiv0\pmod2$.
Then for all
$\,z\in\A(\t^\perp)$
\[ D_c(\t^4\,z)=-4\,D_c(\t^2 \, x\,z)-4\,D_c(\,z). \]
\label{4ofem}\end{cor}
\begin{pf} In $\hat{X}=X\#\CPb$ the class $\t+e$ is represented by a $2$-sphere
of self-intersection
$-3$, and $(\t-2e)\cdot(\t+e)=0$. From Lemma~\ref{blowuplow} we get
\[ \hat{D}_c((\t-2e)^3 \,(\t+e)\,z) = D_c(\t^4\,z)-8\,\hat{D}_c(e^4\,z) =
D_c(\t^4\,z)+16\,D_c(\,z).
\] On the other hand, by Theorems~\ref{RuberSO3} and \ref{3curve} and by
Lemma~\ref{blowuplow},
\begin{eqnarray*}
\hat{D}_c((\t-2e)^3 \,(\t+e)\,z)&=&-\hat{D}_{c+\t+e}((\t-2e)^3\,z)
   =6\,\hat{D}_{c+\t+e}(\t^2\, e\,z)+8\,\hat{D}_{c+\t+e}(e^3\,z)\\
  &=&6\,D_{c+\t}(\t^2\,z)-8\,D_{c+\t}(x\,z) = 12\,D_c(\,z)-4\,D_c(\t^2\, x\,z)
\end{eqnarray*}
and the result follows.
\end{pf}

\bigskip

\section{The blowup equation}

Let $X$ be a simply connected oriented $4$-manifold and let $\hat{X}=X\#\CPb$.
Let $c\in H_2(X;\Z)$. Of course $H_2(\hat{X};\Z)=H_2(X;\Z)\oplus\Z e$ with $e$
the exceptional class.
It follows from Lemma~\ref{blowuplow}(1),(5) that we can write
\[\hat{D}_c=\sum\b_{c,k}\,E^{2k}\]
where $E$ denotes the $1$-form given by $E(y)=e\cdot y$ and
$\b_{c,k}({\a^d})=\hat{D}_c(\a^d \, e^{2k})$ for any $\a\in H_2(X)$. Similarly
\[\hat{D}_{c+e}=\sum\gamma_{c,k}\,E^{2k+1}.\]

Consider $\bar{X}=X\#2\CPb$ with exceptional classes $e_1,\, e_2\in
H_2(\bar{X};\Z)$, and let
$\bar{D}$ denote its Donaldson invariant. Then $e_1+e_2$ has self-intersection
$-2$ and is
represented by an embedded $2$-sphere. Furthermore, the intersection
$(e_1-e_2)\cdot (e_1+e_2)=0$;
so we can apply Corollary~\ref{4ofem} to get
\begin{equation}
\bar{D}_c((e_1-e_2)^r \, (e_1+e_2)^4\,z)=
 -4\,\bar{D}_c((e_1-e_2)^r \, (e_1+e_2)^2\,x\,z)-4\,\bar{D}_c((e_1-e_2)^r\,z)
\label{blowuprecursion}
\end{equation}
for all  $\,z\in \A(X)$.

\begin{lem} There are polynomials, $B_{k}(x)$, independent of $X$, so that for
any
$c\in H_2(X;\Z)$ and $\,z\in \A(X)$ we have $\hat{D}_c(e^k\,z)=
D_c(B_k(x)\,z)$.\end{lem}
\begin{pf} Lemma~\ref{blowuplow}(1) implies that $\b_{c,0}=D_c$. Thus we have
$B_0=1$.
Assume inductively that for $j\le k$, $\hat{D}_c(e^{j}z)= D_c(B_{j}(x)z)$.
Expanding
\eqref{blowuprecursion} via the induction hypothesis we have
\begin{multline*}
\bar{D}_c((e_1-e_2)^{k-3}\,(e_1+e_2)^4\,z)  =
-4\,D_c(\,z\,\sum_{i=0}^{k-3}\binom{k-3}{i}\{ xB_{i+2}(x)B_{k-3-i}(x) \\
-2x B_{i+1}(x)B_{k-2-i}(x)
   + x B_i(x)B_{k-1}(x) + B_i(x)B_{k-3-i}(x)\}\,) = D_c(P(x)\,z)
\end{multline*}
for some polynomial $P$.
On the other hand, expanding the argument of
\[ \bar{D}_c((e_1-e_2)^{k-3}(e_1+e_2)^4\,z) \]
and using the induction hypothesis in a similar fashion, we get
\[\bar{D}_c((e_1-e_2)^{k-3}(e_1+e_2)^4\,z)  = 2\,\hat{D}_c(e^{k+1}\,z) +
D_c(R(x)\,z) \]
for another polynomial $R$. The lemma follows.\end{pf}

\begin{lem} There are polynomials, $S_{k}(x)$, independent of $X$, so that for
any
$c\in H_2(X;\Z)$ and $\,z\in \A(X)$ we have $\hat{D}_{c+e}(e^k\,z)=
D_c(S_k(x)\,z)$.
\end{lem}
\begin{pf} By Theorem~\ref{Ruber} we have for any even $k>0$,
\[
\bar{D}_c((e_1+e_2)^k(e_1-e_2)^2)=2\,D_{c+e_1-e_2}((e_1+e_2)^k)=-2\,D_{c+e_1+e_2}((e_1+e_2)^k) \]
This formula can then be used as above to inductively calculate $S_{k-1}(x)$ in
terms of
$S_1(x),\dots,S_{k-3}(x)$ and $B_0(x),\dots,B_{k+2}(x)$.
\end{pf}

We now explicitly determine the polynomials $B_k(x)$ and $S_k(x)$.
Set \[ B(x,t)=\sum_{t=0}^\infty
B_k(x)\frac{t^k}{k!}\hspace{.25in}\text{and}\hspace{.25in}
   S(x,t)=\sum_{t=0}^\infty S_k(x)\frac{t^k}{k!}\]
Note that
\begin{equation}
\begin{split}\frac{d^n}{dt^n}\,\hat{D}(\exp(te)\,z) & =\hat{D}(\,z\,\sum
e^{k+n}\frac{t^k}{k!})=
D(B_{k+n}(x)\frac{t^k}{k!}\,z) \\ &=\frac{d^n}{dt^n}\,D(B(x,t)\,z) =
D(B^{(n)}(x,t)\,z)
\end{split}
\label{diff}
\end{equation}
where the last differentiation is with respect to $t$.
On $\bar{X} = X\#2\CPb$, we get
$\bar{D}(\exp(t_1e_1+t_2e_2)z)=D(B(x,t_1)\,B(x,t_2)z)$.
Now apply Corollary \ref{4ofem} to $e_1-e_2\in H_2(\bar{X};\Z)$. Since for any
$t\in\R$ the class
$te_1+te_2\in\la e_1-e_2\ra^\perp$, we have the equation
\begin{equation}
\begin{split}\bar{D}(\exp(te_1+te_2)\,(e_1-e_2)^4\,z) & + 4\,
\bar{D}(x\,\exp(te_1+te_2)\,(e_1-e_2)^2\,z)
    \\ &+4\,\bar{D}(\exp(te_1+te_2)\,z) = 0
\end{split}
\label{pre}
\end{equation}
But, for example,
\[ e_1^4\,\exp(te_1+te_2) = \left(\sum e_1^{k+4}\frac{t^k}{k!}\right)
\left(\sum e_2^k\frac{t^k}{k!}\right) =
\frac{d^4}{dt^4}(\exp(te_1))\,\exp(te_2) \]
Arguing similarly and using \eqref{diff}
\begin{equation*}
\begin{split}
& \bar{D}(\exp(te_1+te_2)\,(e_1-e_2)^4\,z)  \\
&=D((2\,B^{(4)}(x,t)\,B(x,t)-8\,B'''(x,t)\,B'(x,t)
   +6\,(B''(x,t))^2)\,z)\\
 &= 2\,D((B^{(4)}\,B-4\,B'''\,B'+3\,(B'')^2 )\,z)
\end{split}
\end{equation*}
where $B=B(x,t)$. Completing the expansion of \eqref{pre} we get
\[2\,D((B^{(4)}\,B-4\,B'''\,B'+3\,(B'')^2 +4x\,(B''B-(B')^2)+2\,B^2)\,z)=0 \]
for all $\,z\in \A(X)$. This means that the expression
\[ B^{(4)}\,B-4\,B'''\,B'+3\,(B'')^2 +4x\,(B''B-(B')^2)+2\,B^2 \]
lies in the kernel of $D:\A(X)\to\R$.

Thus the ``blowup function'' $B(x,t)$ satisfies the differential equation
\[ B^{(4)}\,B-4\,B'''\,B'+3\,(B'')^2 +4x\,(B''B-(B')^2)+2\,B^2 =0 \]
modulo the kernel of $D$.  Of course, the fact that this equation
holds only modulo the kernel of $D$ is really no constraint, since our interest
in
$B(x,t)$ comes from the equation $\hat{D}(\exp(te)z)=D(B(x,t)z)$.

Now let $B=\exp(f(t))$.

\begin{prop} Modulo the kernel of $D$, the logarithm $f(t)$ of $B(x,t)$
satisfies the differential
equation \[ f^{(4)}+6\,(f'')^2+4xf''+2=0\]
with the initial conditions $f=f'=f''=f'''=0$. \ \ \qed \label{DE}\end{prop}

The initial conditions follow from Lemma~\ref{blowuplow}.

\bigskip

\section{The blowup formula}

In order to solve the differential equation of Proposition~\ref{DE}, we set
$u=f''$. Then
 the differential equation becomes
\[ u''+6u^2+4xu+2=0 \] with initial conditions $u(0)=u'(0)=0$.
This is equivalent to the equation
\begin{equation}
(u')^2=-4u^3-4xu^2-4u
\label{uDE}\end{equation}
as can be seen by differentiating both sides of the last equation with respect
to $t$.
Replacing $u$ by $-v$ and completing the cube yields
\[ (v')^2=4(v-\frac{x}{3})^3-\frac43 vx^2 +\frac{4x^3}{27}+4v \]
Finally, letting $y=v-\frac{x}{3}$ we get
\begin{equation}
(y')^2=4y^3-g_2y-g_3 \hspace{.2in} \text{where} \hspace{.2in}
g_2=4\,(\frac{x^2}{3}-1)
\hspace{.1in} \text{and} \hspace{.1in}  g_3={8x^3-36x\over 27} .
\label{yDE}\end{equation}
This is the equation which defines the Weierstrass $\wp$-function. In fact, if
we rewrite \eqref{yDE}
as
\[ \frac{dt}{dy}={1\over\sqrt{4y^3-g_2y-g_3}}\]
then
\[t=\int_y^\infty{ds\over\sqrt{4s^3-g_2s-g_3}}=\wp^{-1}(y)\]
and we see that for arbitrary constants $c$, $y=\wp(t+c)$ gives all solutions
to \eqref{yDE}, and
so $f''=u=-(\wp(t+c)+{x\over3})$ is the general solution of \eqref{uDE}.

The roots of the cubic equation
\[ 4s^3-g_2s-g_3=0 \]
are
\begin{equation}
 e_1=\frac{x}{6}+{\sqrt{x^2-4}\over2}, \hspace{.15in}
e_2=\frac{x}{6}-{\sqrt{x^2-4}\over2},
  \hspace{.15in} e_3=-\frac{x}{3} \label{roots}
\end{equation}
where we have followed standard notation (cf. \cite{Ak}). These correspond to
the half-periods
$\o_i=\wp^{-1}(e_i)$ of the $\wp$-function. The initial condition $f''(0)=0$
implies
that $\wp(c)=-\frac{x}{3}=e_3$; so $c=\o_3+2\varpi$, where
$2\varpi=2m_1\o_1+2m_3\o_3$, with $m_1$,
$m_3\in\Z$, is an arbitrary period. (Note that the initial condition
$f'''(0)=0$ follows
because the half-periods are \,zeros of $\wp'$.) The Weierstrass
\,zeta-function satisfies
$\z'=-\wp$;  thus $f'(t)=\z(t+\o_3+2\varpi)-{tx\over3}+ a$. The constant $a$ is
determined by the
initial condition   $f'(0)=0$;\; $a=-\z(\o_3+2\varpi)$. Since the logarithmic
derivative of the
Weierstrass sigma-function is $\z$, integrating one more time gives
$f(t)=\log\s(t+\o_3+2\varpi)-t\z(\o_3+2\varpi)-{t^2x\over6}+b$,
and the initial condition $f(0)=0$ shows that $b=-\log\s(\o_3+2\varpi)$.
Thus
\[
B(x,t)=e^{f(t)}=e^{-{t^2x\over6}}e^{-t\z(\o_3+2\varpi)}\,{\s(t+\o_3+2\varpi)\over\s(\o_3+2\varpi)}.
\]
For $\o=\o_1$ or $\o_3$ and $\eta=\z(\o)$ we have the formulas
\[\z(u+2m\o)=2m\eta+\z(u)\, ,\hspace{.2in}
\s(u+2m\o)=(-1)^me^{2\eta(mu+m^2\o)}\s(u)\]
(which follow easily from \cite[p.199]{Ak}). Using them, our formula for
$B(x,t)$ becomes
\[ B(x,t)=e^{-{t^2x\over6}}e^{-\eta_3 t}\,{\s(t+\o_3)\over\s(\o_3)}.\]
The above addition formula for the sigma-function implies that
\[\s(t+\o_3)=\s((t-\o_3)+2\o_3)=-e^{2\eta_3t}\s(t-\o_3).\]
Thus
\[ B(x,t)=-e^{-{t^2x\over6}}e^{\eta_3
t}\,{\s(t-\o_3)\over\s(\o_3)}=e^{-{t^2x\over6}}\s_3(t), \]
the last equality by the definition of the quasi-periodic function $\s_3$. In
conclusion,

\begin{thm} Modulo the kernel of $D$, the blowup function $B(x,t)$ is given by
the formula
\[ B(x,t)=e^{-{t^2x\over6}}\s_3(t).\ \ \ \qed\] \label{blowup}\end{thm}

\noindent The indexing of the Weierstrass functions $\s_i$ depends on the
ordering of the roots
$e_i$ of the equation $4s^3-g_2s-g_3=0$. This can be confusing. The important
point is that the
sigma-function we are using corresponds to the root $-{x\over3}$.
One can now obtain the individual blowup polynomials from the formula for
$B(x,t)$. For example,
$B_{12}=-512\,x^4-960\,x^2-408$ and (for fun),
\begin{multline*}
\!B_{30}(x)\!=\!134,217,728\,x^{13}+4,630,511,616\,x^{11}+
  68,167,925,760\,x^9-34,608,135,536,640\,x^7\\
 -39,641,047,695,360\,x^5-9,886,101,110,784\,x^3+543,185,367,552\,x
\end{multline*}
(We thank Alex Selby for help with some computer calculations.)
We also have

\begin{thm} Modulo the kernel of $D$, the blowup function $S(x,t)$ is given by
the formula
\[ S(x,t)=e^{-{t^2x\over6}}\s(t).\] \label{SO3blowup}\end{thm}
\begin{pf} As usual we let $\bar{X}=X\#2\CPb$ with exceptional classes $\ve_1$
and
$\ve_2$. (We have temporarily changed notation to avoid confusion with the
roots $e_i$ of
$4s^3-g_2s-g_3=0$.) Consider the class $\ve_1-\ve_2$ which is
represented by a sphere of self-intersection $-2$. By Theorem~\ref{Ruber} we
have
\[\bar{D}(\exp(t\ve_1+t\ve_2)\,(\ve_1-\ve_2)^2\,z)=2\,D_{\ve_1-\ve_2}(\exp(t\ve_1+t\ve_2)\,z)=
-2\,D_{\ve_1+\ve_2}(\exp(t\ve_1+t\ve_2)\,z)\] for all $\,z\in \A(X)$.
Equivalently we get $D((2\,B''B-2\,(B')^2)\,z)=-2\,D(S^2\,z)$. In other words,
\[ S^2=e^{-{t^2x\over3}}({x\over3}\s_3^2+(\s_3')^2-\s_3\s_3'')\]
Write $\s_3(t)=\exp(h(t))$. Then
\[ S^2= e^{-{t^2x\over3}}e^{2h}({x\over3}-h'') \]
i.e.
\[ S=\pm e^{-{t^2x\over6}}\s_3(t)({x\over3}-h'')^{\frac12}. \]
Since $\exp(h)=\s_3(t)=\s(t)(\wp(t)-e_3)^{\frac12}$, it follows that
$h=\log\s(t)+\frac12\log(\wp(t)-e_3)$. Then
\[h'=\z(t)+\frac12{\wp'(t)\over\wp(t)-e_3}=\z(t)+\frac12(\z(t+\o_3)+\z(t-\o_3)-2\z(t))\]
by \cite[p.41]{Ak}. Thus $h'= \frac12(\z(t+\o_3)+\z(t-\o_3))$, and
\[ h''=\frac12(-\wp(t+\o_3)-\wp(t-\o_3))=-\wp(t+\o_3).\] Thus
\[ S=\pm e^{-{t^2x\over6}}\s_3(t)(\wp(t+\o_3)-e_3)^{\frac12}=
  \pm e^{-{t^2x\over6}}\s_3(t)\left({(e_3-e_1)(e_3-e_2)\over
\wp(t)-e_3}\right)^{\frac12} \]
(\cite[p.200]{Ak}). However, $(e_3-e_1)(e_3-e_2)=1$ (see\eqref{roots}); so
\[ S =\pm e^{-{t^2x\over6}}{\s_3(t)\over\sqrt{\wp(t)-e_3}}=\pm
e^{-{t^2x\over6}}\s(t)\, . \]
To determine the sign, note that (fixing $x$) $S'(0)=S_1=1$ by
Theorem~\ref{blowuplow}(6).
But from our formula $S'(0)=\pm\s'(0)$, and $\s'(0)=1$.
\end{pf}

\bigskip

\section{The blowup formula for manifolds of simple type}

A $4$-manifold is said to be of {\em simple type} \cite{KM} if for all $z\in
\A(X)$ the relation
$D(x^2z)=4\,D(z)$ is  satisfied by its Donaldson invariant. It is clear that if
$X$ has simple type,
then $\hat{X}=X\#\CPb$ does as well. In this case, following \cite{KM}, one
considers the invariant
$\D$ defined by
\[\D(\a)=D((1+{x\over2})\exp(\a))\]
for all $\a\in H_2(X)$. $\D$ is called the {\em Donaldson series} of $X$.
Note that the simple type condition implies that for any $z\in\A(X)$,
\[D((1+{x\over2})\,z\,x)=2\,D((1+{x\over2})\,z)\, ,\]
i.e. $x$ acts as multiplication by $2$ on $D(1+{x\over2})$.
The blowup formula in this case has been determined previously by Kronheimer
and Mrowka.
In this section, we derive that formula by setting $x=2$ in
Theorems~\ref{blowup} and \ref{SO3blowup}. This gives a degenerate case of the
associated
Weierstrass functions. All the formulas below involving elliptic functions can
be found in
\cite{Ak}. The squares $k^2$, ${k'}^2$ of the modulus and complementary modulus
of our Weierstrass
functions are given by
\[ k^2={x-\sqrt{x^2-4}\over x+\sqrt{x^2-4}}  \hspace{.35in}
{k'}^2={2\sqrt{x^2-4}\over x+\sqrt{x^2-4}}.  \]
Thus $k^2=1$ and ${k'}^2= 0$ when $x=2$. The corresponding complete elliptic
integrals of
the first kind are
\begin{eqnarray*}
 K&=&\int_0^1{ds\over\sqrt{(1-s^2)(1-k^2s^2)}}= \int_0^1{ds\over1-s^2}\\
 K'&=& \int_0^1{ds\over\sqrt{(1-s^2)(1-{k'}^2s^2)}}=
\int_0^1{ds\over\sqrt{1-s^2}}
\end{eqnarray*}
Thus $K=\infty$ and $K'={\pi\over2}$ when $x= 2$. Also, when $x=2$ we have
$g_2=\frac43$
and $g_3=-\frac{8}{27}$; so the roots of $4s^3-g_2s-g_3=0$ are
$e_1=e_2=\frac13$ and $e_3=-\frac23$.
This means that when $x=2$ the basic periods are
\[ \o_1={K\over\sqrt{e_1-e_3}}= K=\infty \hspace{.35in}
  \o_3={iK'\over\sqrt{e_1-e_3}}= iK' ={i\pi\over2}\,.\]
In this situation,
\[\s(t)={2\o_3\over\pi}e^{{1\over6}({\pi t\over2\o_3})^2}\sin{\pi t\over2\o_3}=
e^{-{t^2\over6}}\sinh t \]
and
\[\wp(t)=-{1\over3}({\pi\over2\o_3})^2+({\pi\over2\o_3})^2{1\over\sin^2({\pi
t\over2\o_3})}
 ={1\over3}+{1\over\sinh^2t} \, .\]
So
\[ \s_3(t)=\s(t)\sqrt{\wp(t)-e_3}=e^{-{t^2\over6}}\sinh
t\sqrt{1+{1\over\sinh^2t}}
  = e^{-{t^2\over6}}\cosh t. \]

\begin{thm} If $X$ has simple type, the Donaldson series of  $\hat{X}=X\#\CPb$
is
\[ \hat{\D}=\D\cdot e^{-{E^2\over2}}\cosh E \]
where $E$ is the form dual to the exceptional class $e$, i.e. $E(\,z)=e\cdot
\,z$ for all $\,z\in H_2(\hat{X})$. Also
\[ \hat{\D}_e=-\D\cdot e^{-{E^2\over2}}\sinh E\,. \]
\end{thm}
\begin{pf} For $\a\in H_2(X)$ we calculate
\begin{eqnarray*}
\hat{\D}(\a+te) &=&
\hat{D}((1+{x\over2})\exp(\a)\exp(te))=D((1+{x\over2})\exp(\a)B(x,t))\\
   &=& D((1+{x\over2})\exp(\a)e^{-{t^2x\over6}} e^{-{t^2\over6}}\cosh t)
\end{eqnarray*}
The simple type condition implies that $D((1+{x\over2})e^{-{t^2x\over6}})=
D((1+{x\over2})e^{-{t^2\over3}})$.
Hence
\[ \hat{\D}(\a+te)=\D(\a)e^{-{t^2\over2}}\cosh t=\D(\a)(e^{-{E^2\over2}}\cosh
E)(te)\]
as desired. The formula for $\hat{\D}_e$ follows similarly since
$\sinh(E)(te)=-\sinh(t)$.
\end{pf}

A $4$-manifold $X$ is said to have $c$-{\em simple type} if, for $c\in
H_2(X;\Z)$,
$D_c(x^2\,z)=4\,D_c(z)$ for all $z\in \A(X)$.
It is shown in \cite{FS}, and also by Kronheimer and Mrowka, that if $X$ has
simple type, then it has
$c$-simple type for all  $c\in H_2(X;\Z)$. As above we have,

\begin{thm}If $X$ has $c$-simple type,
\begin{eqnarray*}
\hat{\D}_c&=&\D_c\cdot e^{-{E^2\over2}}\cosh E\,. \\
\hat{\D}_{c+e}&=&-\D_c\cdot e^{-{E^2\over2}}\sinh E\,.\ \ \ \qed
\end{eqnarray*}
\end{thm}

 \end{document}